\documentclass[preprint,12pt,3p]{elsarticle}

\usepackage{amssymb}
\usepackage{amsthm,natbib}

\usepackage[utf8]{inputenc}
\usepackage[english]{babel}
\usepackage{amssymb,amsmath,amsthm, mathtools}
\newtheorem{theorem}{Theorem}

\usepackage{appendix}
\usepackage{color, ulem}
\usepackage{pdfsync}
\usepackage{url}
\usepackage{subfigure}
\usepackage{bbold}

\usepackage[utf8]{inputenc}
\usepackage[english]{babel}
\usepackage{amssymb,amsmath, mathtools}

\usepackage{color}
\usepackage{pdfsync}
\usepackage{url}
\usepackage{subfigure}
\usepackage{bbold}
\usepackage{graphicx}
\usepackage{epstopdf}
\newtheorem{assume}[theorem]{Assumption}

\newtheorem{remark}[theorem]{Remark}
\newcommand{\veps}{\varepsilon}

\newcommand{\beq}{\begin{eqnarray*}}
\newcommand{\feq}{\end{eqnarray*}}
\newcommand{\beqn}{\begin{eqnarray}}
\newcommand{\feqn}{\end{eqnarray}}
\newcommand{\bec}{\begin{claim}}
\newcommand{\fec}{\end{claim}}
\newcommand{\becn}{\begin{claim*}}
\newcommand{\fecn}{\end{claim*}}

\journal{Journal }

\begin{document}

\begin{frontmatter}

\title{A Modified Multiple Shooting Algorithm for Parameter Estimation in ODEs Using Adjoint Sensitivity Analysis }
%\tnotetext[label0]{This is only an example}

\author[label1]{Ozgur Aydogmus}
\address[label1]{Social Sciences University of Ankara, Department of Economics\\   Ankara, Turkey }
\ead{ozgur.aydogmus@asbu.edu.tr}
\author[label2]{Ali Hakan TOR}
\address[label2]{Abdullah G\"{u}l University, Department of Applied Mathematics\\  Kayseri, Turkey }
\ead{hakantor@gmail.com}
%\address[label2]{Address Two\fnref{label4}}

%\cortext[cor1]{I am corresponding author}
%\fntext[label3]{I also want to inform about\ldots}
%\fntext[label4]{Small city}

%\ead[url]{author-one-homepage.com}

%\author[label5]{Author Two}
%\address[label5]{Iowa State University, Department Of Mathematics Ames, IA 50014}
%\ead{author.two@mail.com}

%\author[label5]{Author Three}
%\ead{author.three@mail.com}

\begin{abstract}
To increase the predictive power of a model, one needs to estimate its unknown parameters.  Almost all parameter estimation techniques in ordinary differential equation models suffer from either a small convergence region or enormous computational cost. The method of multiple shooting, on the other hand, takes its place in between these two extremes. The computational cost of the algorithm is mostly due to the calculation of directional derivatives of objective and constraint functions. Here we modify { the} multiple shooting algorithm to use the adjoint method in calculating these derivatives. In the literature, this method is known to be a more stable and computationally efficient way of computing gradients of scalar functions. A predator-prey system is used to show the performance of the method and supply all necessary information for a successful and efficient implementation.
\end{abstract}

\begin{keyword} Parameter estimation \sep multiple shooting algorithm \sep adjoint method
%% MSC codes here, in the form: \MSC code \sep code
%% or \MSC[2008] code \sep code (2000 is the default)
\end{keyword}

\end{frontmatter}

%%
%% Start line numbering here if you want
%%
% \linenumbers

%% main text
\section{Introduction}

Modeling by ordinary differential equations (ODEs) is used to accurately depict the physical state of a system in many areas of applied sciences and engineering.  Accurately describing the system and finding its parameters allow for future behavior to be predicted. To estimate parameters in ODEs, we use partially observed noisy data. {Fitting (partially) observed noisy data is a useful way for estimating the unknown parameters  of a system of ODEs.} Implementing such a procedure requires using a convenient ODE solver and an optimization routine.

There are several stochastic and deterministic optimization routines. A detailed discussion of stochastic methods for parameter estimation in ODEs is given in \cite{banga}. Deterministic optimization procedures such as sequential quadratic programming (SQP), Newton methods or quasi-Newton methods can also be employed to estimate the unknown parameters of an ODE system. Contrary to stochastic algorithms, deterministic ones are computationally efficient but tend to converge to local minima.

For parameter estimation problems, the convergence to local minima is prevalent when the single shooting method (also known as initial value approach) is employed. This parameter estimation method uses a single initial condition to produce a trajectory that attempts to fit the noisy data points by minimizing a maximum-likelihood functional with respect to parameters. Hence, the method is computationally efficient but its convergence is highly sensitive to the initial guesses.   

{ The a}bove discussion suggests that there is a trade-off between computational cost and stability of the algorithm. Multiple shooting approach for estimating parameters takes its place in between these two extremes. The method was introduced in \cite{intbook} in 1970s and was substantially enhanced and mathematically analysed in \cite{bock81,bock83}. In this approach, one allows multiple initial values and the problem is considered as a multi-point boundary problem. This is to say that the parameter space is enlarged and discontinuous trajectories are allowed during the optimization process. Employing this method increases the convergence significantly; hence { the} above-mentioned computationally efficient deterministic optimization routines can be used to estimate the parameters. We provide the details of this algorithm in the following section.

We would like to note that { the} multiple shooting method is effective when a gradient based nonlinear programming (NLP) solver such as SQP \cite{sqp,bock4}, quasi-Newton \cite{peifer07} or generalized Gauss-Newton methods \cite{bock3,bock4} is employed. The reason behind the choice of gradient based NLP solvers is that one can compute the sensitivity equations which is useful in finding the gradient of the objective and/or constraint functions.  With that being said, the computational cost of shooting methods is strongly dependent on the efficiency of ODE and sensitivity solvers \cite{booknonlinear}.  

It was claimed in \cite{peifer} that there are three feasible approaches for calculating the derivatives of trajectories with respect to the parameters. These methods are external differentiation, internal differentiation and the simultaneous solutions of the sensitivity equations.  The latter approach is claimed to be the most effective approach among these three. This method requires the simultaneous solution of { the} original ODE system with the sensitivity systems obtained by differentiating the original system with respect to each parameter. If the number of parameters is relatively small this approach may be efficient. On the other hand, the forward sensitivity approach is intractable when the number of state variables and the number of parameters are large. In \cite{cao2,adjoint}, it was claimed that the sensitivities with respect to parameters can be computed more efficiently by employing the adjoint method if the number of parameters is large. Moreover, it was discussed in \cite{efficient} that the adjoint method is the most suitable method to compute sensitivities of a function provided that this function is scalar.

The adjoint method is used in \cite{thesis} to calculate sensitivities for estimating the parameters via single shooting algorithm. In addition, the adjoint and multiple shooting strategies for optimal control of differential algebraic equations systems are combined in \cite{jeon}. To the best of our knowledge, the adjoint method has not been employed to implement a multiple shooting algorithm for estimating the parameters of a system of ODEs. In this study, we aim to modify the classical multiple shooting algorithm so that the adjoint method can be applied efficiently.

The structure of this article is as follows: In the following section, we give a detailed description of the multiple shooting parameter estimation method for systems of ODEs. In Section \ref{sadjmethod}, we modify the multiple shooting algorithm and find the sensitivities with respect to parameters using the adjoint method.  In Section \ref{sexample}, we consider a Lotka-Volterra system and implement the modified algorithm to estimate its parameters from disturbed data. Lastly, we discuss and summarize our findings in Section \ref{sconc}.

\section{The Classical Multiple Shooting Method for Parameter Estimation}\label{sec2}
In this section, we recall { the} classical multiple shooting method whose detailed mathematical analysis was performed in \cite{bock81,bock83}. The method is used to estimate the parameters of the system and some applications of this method to measured data are given in \cite{richter,timmer,stirbet,von}.

We start by considering { a} $d-$dimensional state variable $\mathbf{x(t)\in\mathbb R^d }$ at time $t\in I=[t_0,t_f]$ of a continuous time ordinary differential equation (ODE) { satisfying} the following initial value problem:
\beqn\label{ode}\mathbf{\dot x}(t)=f\big(\mathbf x(t),t,\mathbf p\big), ~~~~\mathbf x (t_0)=\mathbf{x_0}.\feqn We would like to note that the right hand side of the former equality depends on the parameter vector $\mathbf p \in\mathbb{R}^m.$ { The multiple shooting method requires using an enlarged parameter space. This ensures that the procedure has more flexibility for searching the parameter space and circumventing local minima. This enlargement is realized by subdividing time interval with multiple shooting nodes as follows:} \beq t_0=\tau_0<\tau_1<\cdots<\tau_K=t_f.\feq  By introducing the discrete trajectory $\{\mathbf s_j:=\mathbf x(\tau_j)| j=0,1,\cdots, K\},$ one can define the extended parameter vector as:  $\mathbf q=(\mathbf s_0,\mathbf s_1,\cdots,\mathbf s_K,\mathbf p).$ In each subinterval $I_j=[\tau_j,\tau_{j+1}]$ for $j=0,1\cdots,K-1,$ we consider {$K$} independent initial value problems {\beqn\label{diffeqsi}\mathbf{\dot x_j}(t)=f\big(\mathbf x_j(t),t,\mathbf p\big), ~~\mathbf x_j (\tau_j)=\mathbf s_j,~~~~~t\in[\tau_j,\tau_{j+1}]\feqn} where ${ \mathbf x_j} (t)={ \mathbf x_j}(t;\tau_j,\mathbf s_j,\mathbf p)\in\mathbb R^d.$

In \cite{peifer07}, it was required that each subinterval $I_j$ contains at least one measurement. We require that the following assumption holds throughout the paper. 

\begin{assume} \label{ass}Observable variables are indexed by $i=1,2,\cdots,obs$ and measurements are collected at times $\tau_j,$ $j\in M \subset\{0,1,2,\cdots,K\}.$\end{assume} 

Here we remark that the number of observable variables are smaller than or equal to $d.$ We also would like to note that {Assumption \ref{ass}} guarantees that all measurement points are also multiple shooting nodes. {As was noted in \cite{pdebook}, these two sets can be taken equal i.e. $M=\{0,1,2,\cdots,K\}.$} This assumption will be used in the next section to simplify the objective function in such a way that it reduces the computational complexity of the algorithm.  For a detailed discussion, see Remark \ref{rem1}. 

We now introduce the classical multiple shooting algorithm in the following lines following \cite{bock4,peifer07}. Let the measurements for a general function of {the state variables $\mathbf x (t)$ be given as follows}: \beqn\label{measurements}\eta_{ij}=g_{ij}\big(\mathbf x (\tau_j),{\mathbf q})+\veps_{ij}\feqn where measurement errors $\veps_{ij}$ are independent, Gaussian with zero mean and variances {$\sigma_{ij}^2$ for $i\in\{1,2,\cdots,obs\} $ and $j\in M.$ Hence $[\eta_{ij}]$ is a $ obs\times|M|$ matrix} where $|.|$ is used to denote the cardinality of a set. Here the least squares objective function for discontinuous trajectories is defined as follows: \beq L(\mathbf q):=\displaystyle\sum_{i=1}^{obs}\sum_{j\in M}\big(l_{ij}(\mathbf q)\big)^2\feq where $l_{ij}(\mathbf q)=\big(\eta_{ij}-g_{ij}\big({ \mathbf x_j}(\tau_j),\mathbf q\big)\big)/\sigma_{ij}.$ Then we need to consider the following constrained nonlinear optimization problem:
\begin{equation}\label{opt1}
\begin{aligned}
& \underset{\mathbf q}{\text{minimize}}
& & L(\mathbf q)\\
& \text{subject to}
& & {G}_j\big({ \mathbf x_j},\mathbf q\big)=0~~\text{for } j=0,1,\cdots,K-1.
\end{aligned}
\end{equation} where the continuity constrains are given by ${G_j\big( { \mathbf x_j},\mathbf q\big)={ \mathbf x_j}(\tau_{j+1};\tau_j,\mathbf s_j,\mathbf p)-\mathbf s_{j+1} \in \mathbb R^d}.$ Note that each $g_j$ is a $1\times d$ vector and does not specify a scalar valued function. We also note that some optional equality or inequality constraints may be added to this optimization problem \eqref{opt1} (see e.g. \cite{peifer07}).   

Non-linear optimization problem \eqref{opt1} can be solved iteratively using the
generalized-quasi-Newton method \cite[pp.24-25]{book}.  An update step $\Delta\mathbf q=(\Delta \mathbf s_0,\Delta\mathbf s_1,\cdots,\Delta\mathbf s_K,\Delta\mathbf p )$ can be calculated by solving the following Linearized Constrained Least Squares Problem:
\begin{equation}\label{linopt}
\begin{aligned}
& \underset{\mathbf q}{\text{minimize}}
& & \big(l_{ij}(\mathbf q^0)+J_{\mathbf q}l_{ij}(\mathbf q^0)\Delta\mathbf q)\big)^2\\
& \text{subject to}
& & {G}_j\big({ \mathbf x_j},\mathbf q^0\big)+J_{\mathbf q} {G}_j\big({ \mathbf x_j},\mathbf q^0\big)\Delta\mathbf q=0~~\text{for } j=0,1,\cdots,K-1.
\end{aligned}
\end{equation}
for some initial guess $\mathbf q^0=\big(\mathbf s_0^0,\mathbf s_1^0,\cdots,\mathbf s_K^0,\mathbf p^0\big).$ Here $J_q$ denotes the Jacobian with respect to the parameters $\mathbf q$ of the corresponding function. Hence $\Delta \mathbf q$ is used to update parameters as follows: $\mathbf q^l=\mathbf q^{l-1}+\Delta \mathbf q.$ We would like to note that quasi-Newton algorithms has been used to implement multiple shooting procedures. More detailed description of this iterative process is given in \cite[pp.12-17]{bock3}.

Note that one needs to calculate jacobians $J_{\mathbf q}{ G}_j\big({ \mathbf x_j},\mathbf q^0\big)$ in { the} above-given algorithm. This calculation requires finding the following derivatives of the trajectory with respect to parameters:\beqn\label{jac1} \frac{\partial { \mathbf x_j}(\tau_{j+1};\tau_j,\mathbf s_j,\mathbf p)}{\partial p_l} \text{ and } \frac{\partial { \mathbf x_j}(\tau_{j+1};\tau_j,\mathbf s_j,\mathbf p)}{\partial s_{jk}}  \feqn for $l=1,2,\cdots,m$ and $k=1,2,\cdots,d.$  To find these derivatives, one needs to solve sensitivity equations that is a flexible and quite efficient approach (compared to internal or external differentiation \cite{peifer}).  The cost of calculating the sensitivities is the simultaneous integration of the sensitivity equations by solving a system of $m+d$ differential equations at each iteration.

\section{The Modified Algorithm and Sensitivity Analysis Using the Adjoint Method}\label{sadjmethod}

Here our aim is to employ the adjoint method \cite{adjoint} so that the computational cost of the algorithm decreases. First, we would like to note that gradients of objective and constraint functions given in \eqref{opt1} can be calculated using the adjoint method. However, the computational cost of calculating the gradients will be larger than the forward sensitivity analysis since the adjoint method is only effective to find gradients of scalar valued functions as noted in \cite{efficient}. This is to say that one needs to solve $d$ adjoint equations for each one of $d$ elements of vector $g_j$ to calculate the gradient of each $g_j.$ This implies that finding a solution to \eqref{opt1} with no modifications requires solving $d^2$ differential equations in each subinterval per iteration. 

{ The} above discussion implies that the adjoint method cannot be effectively used to calculate the derivatives given in \eqref{jac1}. To reduce the computational cost of the algorithm, we need to consider an equivalent optimization problem. First, recall that one needs to solve a system of $d$ linear differential equations  to find the directional derivatives for each scalar function of the state variables $\mathbf x.$  Hence we need to reduce the number of appearances of the state vectors $\mathbf x$ in objective and constraint functions. 

Using Assumption \ref{ass} along with the initial conditions given in \eqref{diffeqsi} and equation \eqref{measurements}, objective function of \eqref{opt1} can be written as follows: \beqn\label{redefnd} l_{ij}(\mathbf q)=\sigma_{ij}^{-1}\big(\eta_{ij}-g_{ij}\big(\mathbf s_j,\mathbf q\big)\big).\feqn 

\begin{remark}\label{rem1}With this simplification, the objective function does depend on the parameter vector $\mathbf q$ but not {on} the state variables $\mathbf x.$ Hence, there is no need to use the adjoint method when calculating the gradient of the objective function $L(\mathbf q).$ \end{remark}

Now we consider the following optimization problem: \begin{equation}\label{opt2}
\begin{aligned}
& \underset{\mathbf q}{\text{minimize}}
& &  L(\mathbf q)\\
& \text{subject to}
& &h^{(j)}({ \mathbf x_j},\mathbf q)= \|{ G}_j\big({ \mathbf x_j},\mathbf q\big)\|_2^2=0~~\text{for } j=0,1,\cdots,K-1.
\end{aligned}
\end{equation} 
One can easily see that this problem is equivalent to \eqref{opt1}. The difference between these two, on the other hand, {is that there is a scalar constraint} for each subinterval in the latter while the former contains vector valued constraints.

Here our aim is to compute the sensitivities \beqn\label{dsens}\frac{d h^{(j)}}{d{\mathbf q}}\big(\mathbf x_j(\tau_{j+1}),\mathbf q\big) \text{  for  } j=0,1,\cdots,K-1\feqn using the adjoint method. This calculation requires solution of $d$ differential equations on each subinterval rather than $m+d$ differential equations as in the case of sensitivity equations.

\subsection*{The Adjoint Method}
Here we consider the initial value problem \eqref{diffeqsi}. For any interval number $j$ between $0$ and $K-1,$ we have the trajectory ${ \mathbf x_j}(t)={ \mathbf x_j}(t;\tau_i,\mathbf s_i,\mathbf p)$ then, we consider the following function: \beq H_j\big(\mathbf q\big)=\int_{\tau_j}^{\tau_{j+1}}h^{(j)}\big({ \mathbf x_j}(t),\mathbf q\big)\,dt\feq to obtain desired sensitivities \eqref{dsens}. Let $\mathbf{\lambda}(t)$ be any vector valued function of dimension $d$ defined for $t\in[\tau_j,\tau_{j+1}].$ 
Now we consider the following augmented function{:}
\beqn \mathcal L(\mathbf q)=H_j(\mathbf q)+\int_{\tau_j}^{\tau_{j+1}}\lambda^T(t)\big(\mathbf{\dot { \mathbf x}}_j(t)-f\big({ \mathbf x_j}(t),t,\mathbf p\big)\,dt.\label{augmobj}\feqn  Here we use integration by parts to obtain \beq \int_{\tau_j}^{\tau_{j+1}}\lambda^T\mathbf{\dot { \mathbf x}}_j\,dt=\lambda^T{ \mathbf x_j}\Big|_{\tau_j}^{\tau_{j+1}}-\int_{\tau_j}^{\tau_{j+1}}\dot\lambda^T{ \mathbf x_j}\,dt.
\feq Using this equality in \eqref{augmobj}, we obtain \beqn \mathcal L(\mathbf q)=\int_{\tau_j}^{\tau_{j+1}}\Big(h^{(j)}\big({ \mathbf x_j}(t),\mathbf q\big)-\dot\lambda^T(t){ \mathbf x_j}(t)-\lambda^T(t)f\big({ \mathbf x_j}(t),t,\mathbf p\big)\Big)\,dt +\lambda^T{ \mathbf x_j}\Big|_{\tau_j}^{\tau_{j+1}}
\label{obj}\feqn

We would like to note that $\frac{d\mathcal L}{d{\mathbf q}}=\frac{d H_j}{d{\mathbf q}}.$ Hence, taking the total derivative of augmented objective function \eqref{obj}, we obtain:
\beqn  \frac{d\mathcal L}{d{\mathbf q}}=\int_{\tau_j}^{\tau_{j+1}}\Big(\big(h^{(j)}_{{ \mathbf x_j}}-\dot\lambda^T(t)-\lambda^T(t)f_{{ \mathbf x_j}}\big) { \mathbf x_j}_{\mathbf q}+h^{(j)}_{\mathbf q}-\lambda^Tf_{\mathbf q}\Big)\,dt+\big(\lambda^T{ \mathbf x_j}\big|_{\tau_j}^{\tau_{j+1}}\big)_{\mathbf q}.
\label{jacobian1}\feqn
To obtain sensitivity of $H_j(\mathbf q)$ with respect to parameters $\mathbf q$, we require that $\lambda(t)$ satisfies the following initial value problem \beqn\label{sens} \dot\lambda^T(t)=h^{(j)}_{{ \mathbf x_j}}-\lambda^T(t)f_{{ \mathbf x_j}},~~~~ \mathbf \lambda^T(\tau_{j+1})=0 .\feqn  Hence, by \eqref{jacobian1}, sensitivity of $H_j$ can be calculated as follows: \beqn  \frac{d H_j}{d{\mathbf q}}=\frac{d\mathcal L}{d{\mathbf q}}=\int_{\tau_j}^{\tau_{j+1}}\Big(h^{(j)}_{\mathbf q}-\lambda^Tf_{\mathbf q}\Big)\,dt-\lambda^T(\tau_j)\big(\mathbf s_j\big)_{\mathbf q}.\label{jacobian2}\feqn
As noted before, we need to calculate the sensitivities of $h^{(j)}$  at time $\tau_{j+1}$ (see eq. \eqref{dsens}) rather than sensitivities of $H_j.$ By Leibnitz integral rule, we have the following equality: \beq \frac{dh^{(j)}}{d\mathbf q}\big({\mathbf x_j}(\tau_{j+1}),\mathbf q\big)=\frac{d}{d\tau_{j+1}}\frac{dH_j}{d\mathbf q}\feq Using \eqref{jacobian2} in the above equality, we obtain:\beqn \frac{dh^{(j)}}{d\mathbf q}\big({ \mathbf x_j}(\tau_{j+1}),\mathbf q\big)=\Big(h^{(j)}_{\mathbf q}-\lambda^Tf_{\mathbf q}\Big)(\tau_{j+1})-\int_{\tau_j}^{\tau_{j+1}}\lambda_{\tau_{j+1}}^Tf_{\mathbf q}\,dt-{ \lambda}^T_{\tau_{j+1}}(\tau_j)(\mathbf s_j)_{\mathbf q}
\label{sensofh}\feqn where $\lambda_{\tau_{j+1}}^T=\frac{\partial \lambda^T}{\partial\tau_{j+1}}.$ By \eqref{sens}, this newly introduced quantity satisfies the following differential equation: \beqn\label{adjeq}\dot\lambda_{\tau_{j+1}}^T=-\lambda_{\tau_{j+1}}^T f_{{ \mathbf x_j}}.
\feqn By taking the total derivative of the equality $\lambda^T(\tau_{j+1})=0$ as in \cite{adjoint}, we obtain $\lambda_{\tau_{j+1}}^T(\tau_{j+1})=-\dot\lambda^T(\tau_{j+1}).$ This implies by \eqref{sens} that $\lambda_{\tau_{j+1}}^T(\tau_{j+1})=-\dot\lambda^T(\tau_{j+1})=-h^{(j)}_{{ \mathbf x_j}}(\tau_{j+1}).$ Hence, \eqref{sensofh} can be written in the following form: \beqn \frac{dh^{(j)}}{d\mathbf q}({ \mathbf x_j},\mathbf q)|_{t=\tau_{j+1}}=h^{(j)}_{\mathbf q}(\tau_{j+1})-\int_{\tau_j}^{\tau_{j+1}}\lambda_{\tau_{j+1}}^T f_{\mathbf q}\,dt-\lambda_{\tau_{j+1}}^T(\tau_j)(\mathbf s_j)_{\mathbf q}.
\label{sesofh}\feqn { The} above-given formula for the calculation of gradients of constraint functions $h^{(j)}$ contains only $2d+m$ non-zero elements. The fact that gradient vector of $h^{(j)}$ is sparse can be used to efficiently implement an algorithm to solve optimization problem \eqref{opt2}. In the following section, we calculate the sensitivities for a Lotka-Volterra system using { the} above equality.

\section{An Application: Classical Lotka-Volterra Predator Prey System }\label{sexample}
To describe an ecological system consisting of one predator and one prey, one can consider the following model of Lotka and Volterra \beqn\label{lotvolt} \dot x^{(1)}&=&-p_1x^{(1)}+p_2x^{(1)}x^{(2)}\\
\dot x^{(2)}&=&p_3x^{(2)}-p_4x^{(1)}x^{(1)}.\nonumber\feqn The measurements at times $\tau_j=j$ for $j=1,2,\cdots,10$ are simulated by numerical integration of \eqref{lotvolt} with  initial conditions $y_1(0)=0.4$ and $y_2(0)=1$ and parameter values $\mathbf p=(1,1,1,1)$ then they are diturbed  by normally distributed pseudo-random noise ($N(0,\sigma^2)$ with $\sigma^2=0.05.$)

The Lotka-Volterra system has been used as a test example for parameter estimation algorithms (see, for instance, \cite{bock3,book,example2,example1}). The solution has singularities at various combinations of parameters values. { We follow \cite[p.6]{bock3} and take $\mathbf p_0=(0.5, 0.5,0.5,-0.2)$} which results in a pole near $t = 3.3,$ where any ODE solver breaks down \cite{book}. Hence, the single shooting approach  with { the} above-given initial guess $\mathbf p_0$ must fail to estimate the parameters $\mathbf p.$

For { the} above-given parameters, denote the measurements for $x^{(1)}$ and $x^{(2)}$ at $\tau_j$ by $d_j^{1}$ and $d_j^{2},$ respectively. Then the constrained nonlinear minimization problem can be written explicitly as follows : \begin{equation}\label{opt11}
\begin{aligned}
& \underset{\mathbf q}{\text{minimize}}
& & \displaystyle\sum_{j=0}^{10}\Big(\big(s_j^{1}-d_j^{1}\big)^2+\big(s_j^{2}-d_j^{2}\big)^2\Big)\\
& \text{subject to}
& & h^{(j)}({\mathbf x_j},\mathbf q)=0~~\text{for } j=0,1,\cdots,{9}.
\end{aligned}
\end{equation} where $h^{(j)}({\mathbf x_j},\mathbf q)=\bigl(\mathbf x^{(1)}_j(\tau_{j+1})-s^1_{j+1}\bigr)^2+ \bigl(\mathbf x^{(2)}_j(\tau_{j+1})-s^1_{j+1}\bigr)^2.$ 

Here our aim is to calculate the gradient of each function in this nonlinear optimization problem. Since the objective function contains only the parameters of the discrete trajectory $\mathbf s_j$ for $j=0,1,\cdots,10$ one can easily find its gradient. On the other hand, computing the gradient of $h^{(j)}$ with respect to parameters $\mathbf q$ requires using formula \eqref{sesofh}. 

It is obvious from the definition of $h^{(j)}$ and \eqref{sesofh} {that} the gradient of each constraint function has a special form. In particular, one needs to calculate the nonzero elements of the gradient vector $\nabla h^{(j)}$ i.e. the derivatives with respect to $\mathbf s_j,$ $\mathbf s_{j+1}$ and $\mathbf p.$ Finding the derivatives with respect to  $\mathbf s_j$ and  $\mathbf p$  requires integration of the following system of nonautonomous linear differential equations from $\tau_{j+1}$ to $\tau_j$: 
\beq 
\begin{bmatrix}
\dot\Lambda_1 (t)& \dot\Lambda_2(t)
\end{bmatrix}=\begin{bmatrix}
\Lambda_1(t) & \Lambda_2(t)
\end{bmatrix}\begin{bmatrix}
-p_1+p_2x^{(2)}_j(t) &p_2x^{(1)}_j(t)\\
-p_4x^{(2)}_j(t)&p_3-p_4x^{(1)}_j(t)
\end{bmatrix}
\feq satisfying $\Lambda_i(\tau_{j+1})=-2\bigl(\mathbf x^{(i)}_j(\tau_{j+1})-s_{j+1}^i\bigr)$ for $i=1,2.$

 Now one can compute the nonzero elements of the gradient $\nabla h^{(j)}$as follows:
 \beqn 
 \frac{d h^{(j)}}{ds_{j+1}^i}&=&-2\bigl(\mathbf x^{(i)}_j(\tau_{j+1})-s_{j+1}^i\bigr) \text{  for  } i=1,2\nonumber\\
 \frac{d h^{(j)}}{ds_{j}^i}&=&\Lambda_i(\tau_j)  \text{  for  } i=1,2\nonumber\\
  \frac{d h^{(j)}}{dp_1}&=&\int_{\tau_j}^{\tau_{j+1}}\Lambda_1(t)x^{(1)}_j(t)\,dt\nonumber\\
    \frac{d h^{(j)}}{dp_2}&=&-\int_{\tau_j}^{\tau_{j+1}}\Lambda_1(t)x^{(1)}_j(t)x^{(2)}_j(t)\,dt\nonumber\\
    \frac{d h^{(j)}}{dp_3}&=&-\int_{\tau_j}^{\tau_{j+1}}\Lambda_2(t)x^{(2)}_j(t)\,dt\nonumber\\
     \frac{d h^{(j)}}{dp_4}&=&\int_{\tau_j}^{\tau_{j+1}}\Lambda_2(t)x^{(1)}_j(t)x^{(2)}_j(t)\,dt\nonumber
 \feqn

\begin{figure}[h!]

\centering
\subfigure[]{
  \includegraphics[width=7cm, height=7cm]{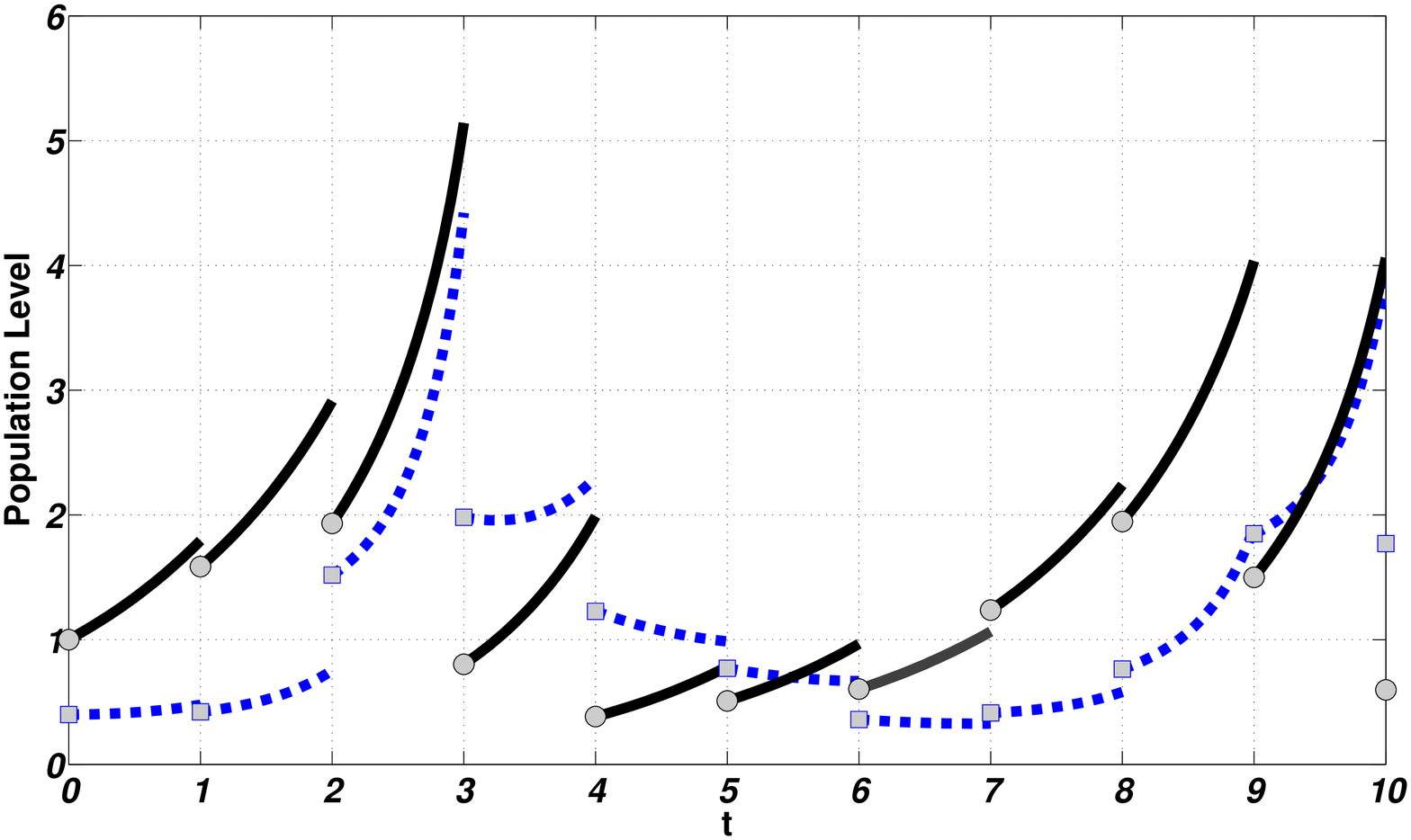}
}
\subfigure[]{
  \includegraphics[width=7cm, height=7cm]{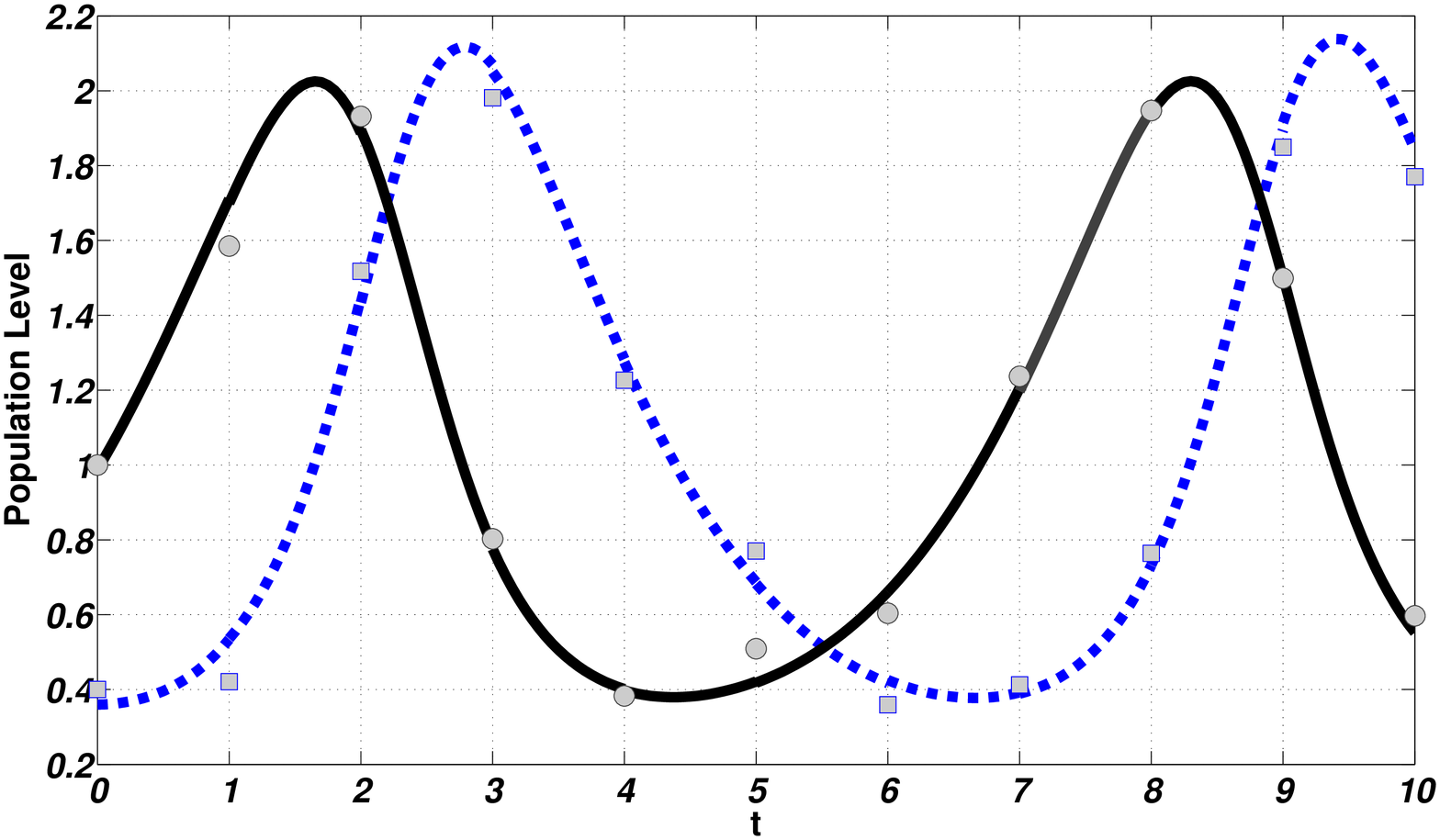}
} 
 
\caption{Multiple shooting trajectories.  Panel (a) illustrates the solutions to \eqref{lotvolt} simulted in each subinterval for the initial guess $\mathbf p_0=(0.5,0.5,0.5,-0.2)$ along with the data points. Similarly, panel (b) shows the same solutions with estimated parameters $p_f=(1.0263, 1.0464 0.9784,0.9688)$ in 8 iterations. }
\label{stationary}
\end{figure}

\begin{remark}
As seen in the above formulas, calculation of directional derivatives with respect to parameters $\mathbf p$ requires integrating some multiplications of state variables $\mathbf x_j(t)$ and $\Lambda_i(t)$ for $i=1,2$ over the interval $I_j.$ Note that using an ODE solver with constant step size provides a computationally efficient way of approximating these integrals.
\end{remark}
We computed the gradients of the constraint functions $h^{(j)}$ using { the} above-defined formulas and obtained the trajectories presented in Figure  \ref{stationary}.  In the panel $(a),$ solutions to the 10 initial value problems with parameters $p_0$ are illustrated as the initial multiple shooting trajectories. In Panel (b), on the other hand, solutions to the same initial value problems are shown for the estimated parameters after 8 iterations of the active set (sqp with a line search) algorithm in Matlab 2015a. Note that this algorithm uses the quasi-Newton approximation of the Hessian matrix; hence it can be used as the optimization routine that is needed for the parameter estimation process. 

We also obtained the following estimations of the parameters $p=(1,1,1,1)$ after 8 iterations of the active set algorithm. We disturbed the data obtained from the simulation of \eqref{lotvolt} with a Gaussian noise term with $\sigma=0.05.$ Repeating parameter estimation process with 10 times with different realizaions of the noise term we obtained the following means and {standard} deviations for the parameters $\mathbf p.$\\
\\
\begin{tabular}{|l|l|l|l|l|}
\hline
              & $p_1$  & $p_2$ & $p_3$  & $p_4$  \\ \hline
mean  & 0.9832 & 0.9870   & 1.0254 & 1.0290 \\ \hline
 standard deviation &0.0554 & 0.0590   & 0.0397  & 0.0371 \\ \hline

\end{tabular} \\
\\
Note that these results are comparable with the results obtained in \cite{bock3}.

%\end{table} 

\section{Conclusion}\label{sconc}

In this paper, we have modified the classical multiple shooting algorithm for parameter estimation in ODEs. In particular, the objective function is taken as a function of only the parameters identifying the discrete trajectory $\mathbf s_j$ and the continuity constraints are taken as scalar functions of the continuous trajectories $\mathbf x_j(t;\tau_j,\mathbf s_j,\mathbf p)$ on each subinterval $I_j.$  This allowed us to use so-called the adjoint method to find the gradient vector of each continuity constraint in each of the intervals $I_j.$ 

Recall that the classical multiple shooting algorithms require solutions to sensitivity equations to determine the Jacobian matrix. Since the number of sensitivity equations increases as the number of parameters increases, the adjoint method provides a computationally efficient alternative for finding such directional derivatives. In addition to this, as claimed in \cite{adjoint}, the adjoint method is more tractable (compared to the forward sensitivity equations) when the number of parameters and state variables {become} large.  Finally, we applied this new method to a well-known predator-prey system which has a singularity for some parameter values. This implies that any ODE solver breaks down, and hence the single shooting method cannot be used to estimate its parameters. We give the explicit expressions for the adjoint equations and nonzero components of the gradient vectors for this system.

Finally, there are several interesting issues that should be further explored or extended.  In this paper, we focused on using the adjoint method to estimate the unknown parameters of a system of ODEs using observed data. We would like to note that both the adjoint method, and { the} multiple shooting method for parameter estimation has been developed for delay differential equations \cite{thesis,calver,voss,voss1}, differential algebraic equations \cite{adjoint,bock6} and partial differential equations \cite{pdeadj,book,multpde,multpde1}. Hence all of these parameter estimation algorithms for qualitatively different systems can be modified, as done in this paper, to find directional derivatives with respect to unknown parameters using the adjoint method.   Moreover, the adjoint method can also be employed to solve constrained optimal control problems (see e.g. \cite{fabien,robot}).

%% References
%%
%% Following citation commands can be used in the body text:
%% Usage of \cite is as follows:
%%   \cite{key}         ==>>  [#]
%%   \cite[chap. 2]{key} ==>> [#, chap. 2]
%%

%% References with bibTeX database:

%\bibliographystyle{elsarticle-num}
% \bibliographystyle{elsarticle-num}
 \bibliographystyle{elsarticle-num-names}
% \bibliographystyle{model1a-num-names}
% \bibliographystyle{model1b-num-names}
% \bibliographystyle{model1c-num-names}
% \bibliographystyle{model1-num-names}
% \bibliographystyle{model2-names}
% \bibliographystyle{model3a-num-names}
% \bibliographystyle{model3-num-names}
% \bibliographystyle{model4-names}

% \bibliographystyle{model6-num-names}

%\bibliography{sample}

\begin{thebibliography}{10}
\expandafter\ifx\csname url\endcsname\relax
  \def\url#1{\texttt{#1}}\fi
\expandafter\ifx\csname urlprefix\endcsname\relax\def\urlprefix{URL }\fi
\expandafter\ifx\csname href\endcsname\relax
  \def\href#1#2{#2} \def\path#1{#1}\fi

\bibitem{banga}
J.~R. Banga, C.~G. Moles, A.~A. Alonso, Global optimization of bioprocesses
  using stochastic and hybrid methods, in: Frontiers in global optimization,
  Springer, 2004, pp. 45--70.

\bibitem{intbook}
J.~Stoer, R.~Bulirsch, Introduction to Numerical Analysis, Springer-Verlag,
  New York, 1993.

\bibitem{bock81}
H.~G. Bock, Numerical treatment of inverse problems in chemical reaction
  kinetics, in: Modelling of chemical reaction systems, Springer, 1981, pp.
  102--125.

\bibitem{bock83}
H.~G. Bock, Recent advances in parameter identification techniques for ODE, in:
  Numerical treatment of inverse problems in differential and integral
  equations, Springer, 1983, pp. 95--121.

\bibitem{sqp}
P.~Drag, K.~Stycze{\'n}, Multiple shooting SQP algorithm for optimal control of
  DAE systems with inconsistent initial conditions, in: Recent Advances in
  Computational Optimization, Springer, 2015, pp. 53--65.

\bibitem{bock4}
H.~G. Bock, S.~K{\"o}rkel, J.~P. Schl{\"o}der, Parameter estimation and optimum
  experimental design for differential equation models, in: Model based
  parameter estimation, Springer, 2013, pp. 1--30.

\bibitem{peifer07}
M.~Peifer, J.~Timmer, Parameter estimation in ordinary differential equations
  for biochemical processes using the method of multiple shooting, IET Systems
  Biology 1~(2) (2007) 78--88.

\bibitem{bock3}
H.~G. Bock, E.~Kostina, J.~P. Schl{\"o}der, Direct multiple shooting and
  generalized Gauss-Newton method for parameter estimation problems in ODE
  models, in: Multiple Shooting and Time Domain Decomposition Methods,
  Springer, 2015, pp. 1--34.

\bibitem{booknonlinear}
L.~T. Biegler, Nonlinear programming: concepts, algorithms, and applications to
  chemical processes, Vol.~10, Siam, 2010.

\bibitem{peifer}
M.~Peifer, J.~Timmer, Parameter estimation in ordinary differential equations
  using the method of multiple shooting—a review (2005).

\bibitem{cao2}
Y.~Cao, S.~Li, L.~Petzold, Adjoint sensitivity analysis for
  differential-algebraic equations: algorithms and software, Journal of
  Computational and Applied Mathematics 149~(1) (2002) 171--191.

\bibitem{adjoint}
Y.~Cao, S.~Li, L.~Petzold, R.~Serban, Adjoint sensitivity analysis for
  differential-algebraic equations: The adjoint DAE system and its numerical
  solution, SIAM Journal on Scientific Computing 24~(3) (2003) 1076--1089.

\bibitem{efficient}
B.~Sengupta, K.~J. Friston, W.~D. Penny, Efficient gradient computation for
  dynamical models, NeuroImage 98 (2014) 521--527.

\bibitem{thesis}
J.~J. Calver, Parameter estimation for systems of ordinary differential
  equations, Ph.D. thesis (2019).

\bibitem{jeon}
M.~Jeon, Parallel optimal control with multiple shooting, constraints
  aggregation and adjoint methods, Journal of Applied Mathematics and Computing
  19~(1-2) (2005) 215.

\bibitem{richter}
O.~Richter, P.~N{\"o}rtersheuser, W.~Pestemer, Non-linear parameter estimation
  in pesticide degradation, Science of the Total Environment 123 (1992)
  435--450.

\bibitem{timmer}
J.~Timmer, H.~Rust, W.~Horbelt, H.~Voss, Parametric, nonparametric and
  parametric modelling of a chaotic circuit time series, Physics Letters A
  274~(3-4) (2000) 123--134.

\bibitem{stirbet}
A.~D. Stirbet, P.~Rosenau, A.~C. Str{\"o}der, R.~J. Strasser, Parameter
  optimisation of fast chlorophyll fluorescence induction model, Mathematics
  and Computers in Simulation 56~(4-5) (2001) 443--450.

\bibitem{von}
H.~von Gr{\"u}nberg, M.~Peifer, J.~Timmer, M.~Kollmann, Variations in
  substitution rate in human and mouse genomes, Physical Review Letters 93~(20)
  (2004) 208102.

\bibitem{pdebook}
S.~B. Hazra, Large-scale PDE-constrained optimization in applications, Vol.~49,
  Springer Science \& Business Media, 2009.

\bibitem{book}
K.~Schittkowski, Numerical data fitting in dynamical systems: a practical
  introduction with applications and software, Vol.~77, Springer Science \&
  Business Media, 2013.

\bibitem{example2}
J.~Swartz, H.~Bremermann, Discussion of parameter estimation in biological
  modelling: Algorithms for estimation and evaluation of the estimates, Journal
  of Mathematical Biology 1~(3) (1975) 241--257.

\bibitem{example1}
L.~Edsberg, P.-{\AA}. Wedin, Numerical tools for parameter estimation in
  ODE-systems, Optimization Methods and Software 6~(3) (1995) 193--217.

\bibitem{calver}
J.~Calver, W.~Enright, Numerical methods for computing sensitivities for ODEs
  and DDEs, Numerical Algorithms 74~(4) (2017) 1101--1117.

\bibitem{voss}
H.~Voss, M.~Peifer, W.~Horbelt, H.~Rust, J.~Timmer, G.~Gousebet, Identification
  of chaotic systems from experimental data, Chaos and its
  reconstruction’(Nova Science Publishers Inc., New York, 2003) (2003)
  245--286.

\bibitem{voss1}
W.~Horbelt, J.~Timmer, H.~U. Voss, Parameter estimation in nonlinear delayed
  feedback systems from noisy data, Physics Letters A 299~(5-6) (2002)
  513--521.

\bibitem{bock6}
H.~Bock, M.~Diehl, D.~Leineweber, J.~Schl{\"o}der, A direct multiple shooting
  method for real-time optimization of nonlinear DAE processes, in: Nonlinear
  model predictive control, Springer, 2000, pp. 245--267.

\bibitem{pdeadj}
A.~M. Bradley, Pde-constrained optimization and the adjoint method, Tech. rep.,
  Technical Report. Stanford University. https://cs. stanford. edu/\~{}
  ambrad~… (2013).

\bibitem{multpde}
X.~Xun, J.~Cao, B.~Mallick, A.~Maity, R.~J. Carroll, Parameter estimation of
  partial differential equation models, Journal of the American Statistical
  Association 108~(503) (2013) 1009--1020.

\bibitem{multpde1}
T.~M{\"u}ller, J.~Timmer, Parameter identification techniques for partial
  differential equations, International Journal of Bifurcation and Chaos
  14~(06) (2004) 2053--2060.

\bibitem{fabien}
B.~C. Fabien, Numerical solution of constrained optimal control problems with
  parameters, Applied Mathematics and Computation 80~(1) (1996) 43--62.

\bibitem{robot}
M.~Diehl, H.~G. Bock, H.~Diedam, P.-B. Wieber, Fast direct multiple shooting
  algorithms for optimal robot control, in: Fast motions in biomechanics and
  robotics, Springer, 2006, pp. 65--93.

\end{thebibliography}
%\bibliographystyle{elsarticle-harv}
%\bibliographystyle{elsarticle-num}

\end{document}